\RequirePackage{fix-cm}

\documentclass[nospthms,natbib]{svjour3}       
\def\cite#1{\citep{#1}}

\makeatletter \providecommand{\@LN}[2]{} \makeatother

\usepackage{bm}
\usepackage{booktabs} \usepackage{siunitx} \usepackage{pgfplotstable} \pgfplotsset{compat=1.17}
\usepackage[cmex10]{amsmath}
\usepackage{amsfonts,amssymb,amsthm,mathtools}
\usepackage{algorithmic}
\usepackage{algorithm}
\usepackage{array}
\usepackage{textcomp}
\usepackage{stfloats}
\usepackage{url}
\usepackage{verbatim}
\usepackage{graphicx}
\usepackage{subfig}
\usepackage{makecell} 
\usepackage[english]{babel}

\usepackage[pdftitle={Optimizing Flexibility in Power Systems by Maximizing the Region of Manageable Uncertainties},pdfauthor={Aron Zingler Stephane Fliscounakis, Patrick Panciatici, Alexander Mitsos},breaklinks]{hyperref}
\hypersetup{hidelinks,
  colorlinks=true,
  allcolors=black,
  pdfstartview=Fit,
  breaklinks=true}
  
\usepackage{cleveref}
\crefformat{equation}{(#2#1#3)}
\Crefformat{equation}{(#2#1#3)}
\crefformat{figure}{Fig.~#2#1#3}
\Crefformat{figure}{Fig.~#2#1#3}
\crefformat{table}{Table~#2#1#3}
\Crefformat{table}{Table~#2#1#3}
\crefformat{section}{Section~#2#1#3}
\Crefformat{section}{Section~#2#1#3}
\usepackage{mymath}

\usepackage{xargs}

\usepackage{etoolbox}

\makeatother

\newcommand{\maxflow}{maximal additional power transfer capacity}
\newcommand{\SVT}{
\begin{footnotesize}
Process Systems Engineering (AVT.SVT), RWTH Aachen University, 52074 Aachen, Germany
\end{footnotesize}\\
}
\newcommand{\JARACSD}{
\begin{footnotesize}
JARA-CSD, 52056 Aachen, Germany
\end{footnotesize}\\
}

\newcommand{\IEK}{
\begin{footnotesize}
Institute of Energy and Climate Research: Energy Systems Engineering (IEK-10), Forschungszentrum J{\"u}lich GmbH, 52425 J{\"u}lich, Germany.
\end{footnotesize}\\
}

\newcommand{\RTE}{
\begin{footnotesize}
R\'eseau de transport d'\'electricit\'e, 
Paris, France.
\end{footnotesize}\\
}
\begin{document}
\title{Optimizing Flexibility in Power Systems by Maximizing the Region of Manageable Uncertainties}
\titlerunning{Optimizing Flexibility in Power Systems}

\newcommand{\orcid}[1]{\href{https://orcid.org/#1}{\textsuperscript{\includegraphics[height =  1.9ex]{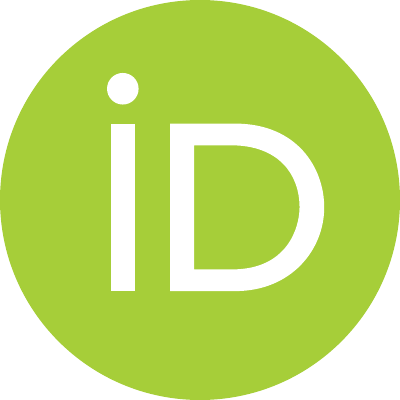}}}}
\author{Aron Zingler \orcid{0000-0003-1135-0236} \and St\'ephane Fliscounakis \and Patrick Panciatici \orcid{0009-0000-8936-094X} \and Alexander Mitsos \orcid{0000-0003-0335-6566}}
\authorrunning{Aron Zingler \and St\'ephane Fliscounakis \and Patrick Panciatici \and Alexander Mitsos}

\institute{
A. Mitsos (corresponding author)\\
 \email{amitsos@alum.mit.edu},
\at \JARACSD
\and A. Mitsos  \at \IEK
\and A.Zingler \and A. Mitsos \at \SVT
\and St\'ephane Fliscounakis \and Patrick Panciatici \at \RTE
              \begin{acknowledgements}
              This research is funded by R\'eseau de transport d'electricit\'e (RTE, France) through the project “Hierarchical Optimization for Worst-Case Analysis of Power Grids”.
              This preprint has not undergone peer review (when applicable) or any post-submission
			  improvements or corrections. The Version of Record of this article is published in Optimization and Engineering, and is available online at https://doi.org/10.1007/s11081-025-09958-z”.
			  \end{acknowledgements}
}

\date{25.11.2024}

\journalname{arxiv.org}

\maketitle

\begin{abstract}
Motivated by the increasing need to hedge against load and generation uncertainty in the operation of power grids, we propose flexibility maximization during operation.
 We consider flexibility explicitly as the amount of uncertainty that can be handled  while still ensuring nominal grid operation in the worst-case.
We apply the proposed flexibility optimization in the context of a DC flow approximation. By using a corresponding parameterization, we can find the maximal range of uncertainty and a range for the manageable power transfer between two parts of a network subject to uncertainty.
 We formulate the corresponding optimization problem as an (existence-constrained) semi-infinite optimization problem and specialize an existing algorithm for its solution. 
 
\end{abstract}

\section{Introduction}
Uncertainty is a primal concern in the operation of power grids since exact values for supply and demand are rarely known in advance.
To ensure robust grid operation, it is therefore important to quantify and improve the flexibility of a grid, i.e., the amount of uncertainty allowed in the grid while maintaining safe operation. 
Although assessing and optimizing the flexibility of power systems is actively explored in the literature for different notions of flexibility \cite{nosair2015flexibility,ma2013evaluating}, we want to investigate flexibility in the specific sense of an explicit measure of the range of manageable uncertain values. Manageable uncertainties are uncertainties for which operational constraints on the grid can be ensured by taking appropriate control actions. For example, a grid operator might be interested in changing generator set-points to maximize the range of power injections that can be handled by the load distribution of the grid without violating the power limit on critical lines.  This perspective on flexibility seems to be underutilized in the power system literature. 
In contrast, this flexibility concept was already introduced in the 80s within process systems engineering \cite{grossmann1983optimization}. However, as pointed out by \cite{wei2012flexibility}, the number of uncertain variables  was typically small in these studies, whereas power systems applications require many variables.

We introduce a method for maximizing flexibility of power systems under uncertainty in a worst-case framework.
We focus on situations where an explicit description of the region of \emph{manageable} uncertainty is not known apriori. We allow for a two-stage approach with mixed integer control decisions. For a given parametrization of the uncertainty region, i.e., the set of all possible uncertainties, we aim at finding preventive actions such that all uncertainties inside that parametrized region are manageable and the region is as large as possible.
In our case studies, we use two parametrizations that allow for an intuitive interpretation of the result.
In the first, we create an inner approximation of the region of manageable uncertainty by searching for the largest scaled hyperbox of uncertainties that can be managed. This allows us to find the largest amount of relative uncertainty of the node injections that can be handled.
In the other, we parametrize the uncertainty region based on the amount of additional power flow from one region to another caused by the uncertainty compared to  a forecasted base case. Here the aim is to find a \maxflow{} where every additional power flow up to the capacity is manageable.
\section{Motivating example \label{sec:motivating_example}}
\begin{figure}
    \centering
    \includegraphics[width=0.6\textwidth]{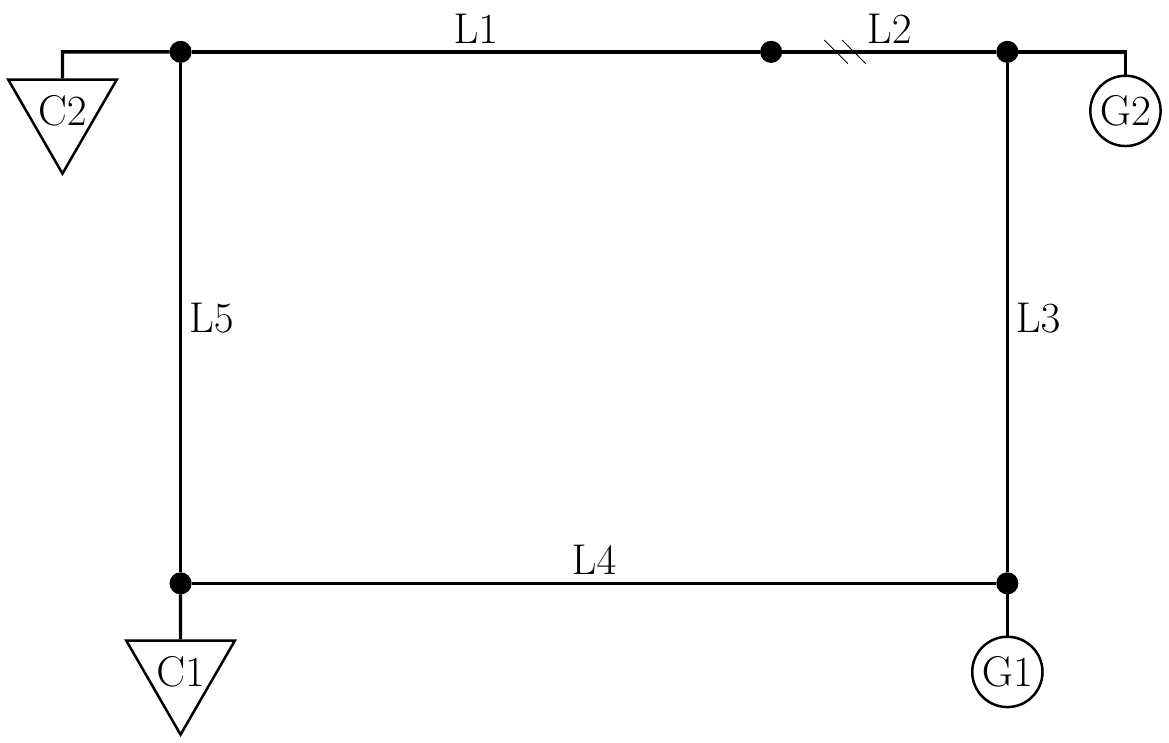}
    \caption{Motivating power grid  instance with two generators and two uncertain loads.}
    \label{fig:min_example_flex_grid}
\end{figure}

Before discussing the general problem formulation and concrete model equations, we give a simple example of a question that can be answered using our approach. 
Consider the power grid instance in \cref{fig:min_example_flex_grid}.

The lines $L_i$, $i=1,2,3,4$,  have equal admittance $h$  and  limit ${P}^{lim}$ for maximal power flow across the line.

The loads $I_{c1}$ and $I_{c2}$ at the two consumers $C1$ and $C2$ have predicted values $I_{c1}^{pred}$ and $I_{c2}^{pred}$, but are uncertain. The loads may vary between the bounds $I_{c1}\in[I_{c1}^{pred}-\delta,I_{c1}^{pred}+\delta]$ and $I_{c2}\in[I_{c2}^{pred}-2\delta,I_{c2}^{pred}+\delta]$.

We want to set the set-points $I_{g1}^0$ and $I_{g2}^0$ of the two generators $G1$ and $G2$ in such a way that the predicted loads are balanced and the limits on the lines are not violated for all considered uncertain load values.
To balance the difference between predicted and encountered loads, we prescribe a load distribution strategy that increases the injections at both generators by the same amount, i.e., $I_{gi}=I_{gi}^0+\Delta I_{gi},\; i=1,2$ and $\Delta I_{g1}=\Delta I_{g2}$.  Additionally, we can choose to open line $L_2$ as a control action.

As an index of flexibility, we are interested in the largest value of $\delta$ that still admits a feasible solution that is robust in the worst-case under the DC flow approximation.

In the following, we describe an approach  to optimize this index, i.e.  to calculate the optimal $I_{gi}^0$ leading to the largest corresponding $\delta$, in this example and in more general cases.

\section{Literature review}
The importance of flexibility in power systems is well recognized in the literature\cite{babatunde2020power,oecd2011harnessing,akrami2019power}. The problem of flexibility has multiple facets \cite{akrami2019power} and the way flexibility is measured  can range from ``megawatts (MW) available for ramping up and down over time''\cite{oecd2011harnessing} over the fraction of the operating region that stays  feasible when taking network limitations into account \cite{gomez2019operational} to the ``sum of probabilities of the scenarios satisfied by the chosen solution strategy'' \cite{menemenlis2011thoughts}.
For the purpose of this work, we define flexibility as an explicit index that describes which values of the uncertain parameters can be handled.

We are not the first to point out the potential of this idea for power systems.
 In \cite{wei2012flexibility}, the authors already introduce this idea of flexibility analysis from process systems to power systems.  
 However, instead of treating uncertainty explicitly, they increase robustness with an approach that amounts to finding the maximal restriction to selected inequality constraints. 
 
Other authors have considered the uncertainty explicitly, but evaluated the flexibility of a system without trying optimize parameters to maximize it:
In \cite{gomez2019operational}, the authors approximate the size of the region of manageable uncertainty, but focus on cases with access to an explicit half-space representation.
The authors of \cite{bucher2015quantification} investigate a similar approach.
They project an explicit half-space representation of the region of feasible controls for given disturbances into the space spanned by manageable disturbances. 
In \cite{zhao2015unified}, the authors pursue a goal similar to ours: they determine the maximal size of a hyperbox containing manageable uncertainties while considering control actions. In contrast to their approach, we consider a general parameterization of the shape of uncertainty region and consider not only the determination, but also optimization of this flexibility metric. 

We aim for rigorous upper and lower bounds of our flexibility metric. As a result, our approach requires a global solution of the resulting subproblems, which becomes computationally expensive for nonlinear models. Furthermore, convergence can only be guaranteed when the grid state is uniquely defined for fixed values of preventive actions and uncertainties. We restrict ourselves to the DC flow approximation for these reasons, similarly to the references above.

The heuristic evaluation of a flexibility metric similar to ours for fixed preventive actions is investigated by the authors of \cite{capitanescu2021power} in the context of nonconvex AC power flow. A heuristic procedure is applied that does not give any guarantees on convergence or on optimality of the result. 
 In \cite{lee2021robust}, the authors consider the nonconvex AC-power equations and construct an inner convex approximation of the feasible region around an initial point under uncertainty. For uncertainty, they consider uncertain power injections $\mb{y}$ that are known to occur in an ellipse around a nominal point $\mb{y}^0$ with a given radius $\gamma$. As one application, they consider finding locally the maximum robustness margin which is the size of the uncertainty set
against which the solution is robust, i.e., the maximum radius $\gamma$. 
This is the only publication known to us that also considers taking actions to optimize the flexibility of the network in a way that allows for an interpretation of the result in terms of manageable uncertainty values.

\section{Problem Formulation and Solution}
In our general problem formulation, we maximize a scalar index for flexibility $\delta$ by adjusting preventive actions for grid operation $\mb{x}\in \mc{X}$. The flexibility index $\delta$ determines the size of the parametric uncertainty region $\mc{T}(\delta,\mb{x})$.
We must guarantee the existence of  controls $\mb{z} \in \mc{Z}(\mb{x},\mb{y})$ that allow safe grid operation, captured by a continuous inequality constraint $g^s$, for all uncertainty values $\mb{y}\in \bar{\mc{Y}}$  in the parametric uncertainty region $\mc{T}(\delta,\mb{x})$, as well as for a predicted state with fixed $\mb{y}^0$ and $\mb{z}^0$.  We assume that the grid states are uniquely given by a continuous function  $\mb{s}(\mb{x},\mb{y},\mb{z})$ and that the set of controls $\mc{Z}(\mb{x},\mb{y})$ is not empty on $\mc{X}\times \bar{\mc{Y}}$.  

The resulting problem  
\begin{equation}
\begin{aligned} \label{eq:flex}
\inf_{\delta\in \R^+,\mb{x}\in \mc{X}} \; & -\delta\\
\text{s.t. } & g^s(\mb{x},\mb{y}^0,\mb{z}^0,\mb{s}(\mb{x},\mb{y}^0,\mb{z}^0)) \le 0, \\
	        & \forall \mb{y} \in \mc{T}(\delta,\mb{x})\left[ \exists {\mb{z}\in \mc{Z}(\mb{x},\mb{y}}) \right. :\\ 
	        & \qquad \left. g^s(\mb{x},\mb{y},\mb{z},\mb{s}(\mb{x},\mb{y},\mb{z})) \le 0 \right]
\end{aligned}
\end{equation}
is a generalized existence-constrained semi-infinite optimization problem. In our case, the model is based on the DC flow approximation (which will lead to MILP subproblems) and introduced next together with our choices for the preventive actions $\mb{x}$, the uncertainties $\mb{y}$ and the controls $\mb{z}$.  In this context, the function $\mb{s}$ is not known explicitly but is rather given implicitly by the model constraints. To solve this challenging problem, we specialize the discretization algorithm from \cite{djelassi2021esipglobal}, which allows for the inclusion of the implicit function $\mb{s}$ \cite{djelassi2019discretization}.

\section{Model equations}
This work builds upon the model described in \cite{djelassi2018hierarchicalgrid,djelassi2020discretization}.  An overview of the model is given below. For a complete description, we refer to the aforementioned references.
\subsection{Notation}
The power grid topology is represented by a directed graph with nodes $n \in \mc{N}$ and edges $e \in \mc{E}$.  Each edge $e$ has an originating node $\tilde{n}(e)$ and a terminating node $\tilde{m}(e)$.
Each node $n$ has a  set of incoming edges $\mc{E}^+(n)$ and outgoing edges $\mc{E}^-(n)$.
Some of the nodes are connected to generators $g\in \mc{G}$. We denote this set of generator nodes $\mc{N}_G \subseteq \mc{N}$ and the mapping of generators to corresponding nodes as
$n_G(g)$.
Some of the edges are associated with a phase-shifting transformer (PST).  The corresponding set of shifter edges is denoted by $\mc{E}_S \subseteq \mc{E}$.

\subsection{Summary of the model}
In the general formulation \cref{eq:flex}, we differentiate between preventive actions $\mb{x}$ and control actions $\mb{z}$.  They differ in the following way:
The preventive actions $\mb{x}$ have to be selected in a first stage without prior knowledge of the uncertainty values $\mb{y}$.  
In our model, the preventive actions are the set points for power generated by the generators $\mb{x}^T=[I_{1}^{x},...,I_{|G|}^{x}]^T$. 

The value of the flexibility index $\delta$ also has to be selected in the first stage. It determines how much we restrict the uncertainty values that can occur in the second stage. Higher values of $\delta$ increase the size of the parametric uncertainty region $\mc{T}(\delta, \mb{x})$ but also improve the value of the objective function.

In the second stage,  worst-case un\-cer\-tain\-ty values $\mb{y}$ are de\-ter\-min\-ed within the para\-metric un\-cer\-tainty re\-gion $\mc{T}(\delta,\mb{x})$.
They  consist of the uncertain derivations for the nodal power injections $\Delta I_n^y$ from their predicted values $I_n^0$ at nodes $n\in \mc{N}$, i.e., $\mb{y}^T=[\Delta I_1^y,...,\Delta I_{|\mc{N}|}^y]^T$. 

The control variables $\mb{z}$ are chosen in the last stage with full knowledge of the uncertainty values. They include the setting of phase shifter $m_e$ on the edges $e\in \mc{E}$ and binary decision variables $p_b$ that decide if designated pairs of busses $b\in \mathcal{B}$ are merged. In other words, $\mb{z}^T=[\Delta \theta_1,..., \Delta \theta_{|\mc{E}|}, p_1,...,p_{|\mc{B}|}]^T$.

We also consider a fixed load distribution scheme among the generators. The power generation at the generators can be changed by control offsets  $ \Delta I_{g}^{z}$ at generators $g\in \mc{G}$. Those offsets are conceptually control decisions. 
However, as detailed later, because  they are uniquely determined as a function of the uncertainty offsets for the nodal power
injections $\Delta I _n^y$ and the generator set points $I_{g}^{x}$, we should calculate them alongside the uncertainty values \cite{djelassi2019discretization}.

The vector of implicit state variables $\mb{s}$ contains the voltage angles at nodes $n\in\mc{N}$, denoted as $\theta_n$, and the  $\mb{s}=[\theta_1,...,\theta_{|\mc{N}|}]$. The implicit state variables are given uniquely by the implicit function $\mb{s}(\mb{x},\mb{y},\mb{z})$, i.e., they are  determined by the other variables, and the grid equations.

\subsection{Components}
\subsubsection{Edges}
\begin{figure}
    \centering
    \includegraphics{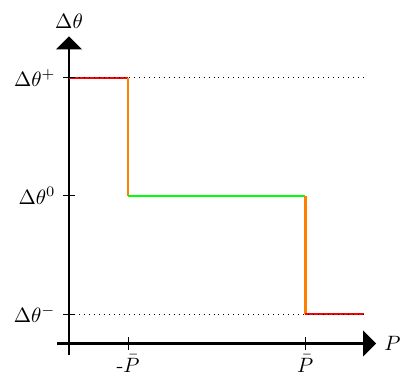}
    \caption{Illustration of the behavior of phase shifters.}
    \label{fig:illustrate_phase_shifter}
\end{figure}

To calculate the power flow $P_e$ over an edge $e\in \mc{E}$, the DC power flow approximation 
\begin{equation}
P_e = h_e(\theta_{\tilde{n}(e)} - \theta_{\tilde{m}(e)} + \Delta\theta_e)
\end{equation}
is used with the positive admittance $h_e$, the voltage angle of the originating node $ \theta_{\tilde{n}(e)}$ and the voltage angle of the terminating node $\theta_{\tilde{m}(e)}$.
The phase shift $\Delta \theta_e$ is zero for all edges $e\in \mc{E}\setminus \mc{E}_S$ without a PST.  For the edges with a PST, the phase-shift of a PST only deviates from zero if  the magnitude of power flow $|P_e|$ reaches or exceeds a specified activation threshold $\bar{P}$. The phase-shift will vary in its bounds, i.e., $[\Delta \theta_e^-,\Delta \theta_e^+]$ to keep the power flow from exceeding the limit until this is prevented by said bounds. This piecewise linear model was introduced in \cite{djelassi2018hierarchicalgrid} to reduce unrealistic continuous control in the PST and is illustrated in \Cref{fig:illustrate_phase_shifter}.

An additional control mechanism is given by bus merging. For a predefined set of pairs of nodes $\mc{B}\subseteq \mc{N} \times \mc{N}$, we allow a single pair to be merged. For each pair $b\in \mc{B}$, the binary variable $p_b$ encodes if the pair of nodes is merged. We can model bus merging with a line that has infinite admittance if closed, but can also be open.
Thus for these lines, we can model their behavior with the following disjunction:

\begin{equation*}
\begin{aligned} 
p_b = 0 &\implies 0 = P_b \\
p_b = 1 &\implies 0 = \theta_{\tilde{n}(b)} - \theta_{\tilde{m}(b)}
\end{aligned}
\end{equation*}

\subsubsection{Nodes}
The injection $I_n$ at a any node $n\in \mc{N}$ must equal the outgoing power-flows over the connected edges

\begin{equation}
    I_n=\sum_{e\in \mathcal{E}: n=\tilde{n}(e)} P_e - \sum_{e\in \mathcal{E}: n=\tilde{m}(e)} P_e .
\end{equation}

The injection $I_n$ at non-generator nodes $n\in \mc{N} \setminus \mc{N}_G$ is uncertain  and given by
\begin{equation}
I_n = I_n^0 + \Delta I _n^y
\end{equation}
with an initial injection of $I_n^0$ and an uncertainty offset $\Delta I _n^y$ within given bounds ($\Delta I _n^y \in [\Delta I _n^{y-},\Delta I _n^{y+}]$).

The injection $I_n$ at nodes $n=n_G(g)$ with generators $g\in \mc{G}$ is given by
\begin{equation}
I_{n_G(g)}= I_{n_G(g)}^0+I_{g}^{x}+\Delta I_{g}^{z}+ \Delta I_{n_G(g)}^{y}
\end{equation}
where $I_{g}^{x}$ is injection by the generator $g$ at its set point and $\Delta I_{g}^{z}$ is the control offset due to load distribution.

The control offset of generators follows a piecewise linear profile. The sum of all injections due to control offsets must cancel the additional injections due to uncertainty.
 In total, the total injection demand
\begin{equation}
\Delta I^t= -\sum_{n \in \mc{N}} \Delta I_{n}^{y} =  \sum_{g \in \mc{G}} \Delta I_{g}^{z}
\end{equation}
must be  injected by load distribution.

Ideally, this total injection would be divided among the generators $g\in \mc{G}$ according to a contribution factor $c_g\ge0$ with $\sum_{g\in \mc{G}} c_g =1$.
However, the demanded control offset $ \Delta I_{g}^{z,dem}=c_g \Delta I^t$ at each generator $g\in \mc{G}$ can be prohibited because the total generator output is limited by the generator bounds $I_{g}^{x-}$ and $I_{g}^{x+}$, i.e., 
\begin{equation}
I_{g}^{x-} \le I_{g}^{x}+\Delta I_{g}^{z} \le I_{g}^{x+} \label{eq:total_gen_limit}.
\end{equation}

In this situation, we introduce an artificially increased total injection demand $\Delta I^{t,inc}$ such that the control offsets resulting from this increased total injection demand and the limits on generator output fulfill
\begin{equation}
\Delta I^t =  \sum_{g \in \mc{G}} \underbrace{ \operatorname{mid}(I_{g}^{x-}, I_{g}^{x}+c_g\Delta I^{t,inc},I_{g}^{x+})  - I_{g}^{x} }_{\Delta I_{g}^{z}}.
\end{equation}
Here, the $\operatorname{mid}$ function returns the median of its elements.

\begin{figure}
    \centering
    \includegraphics{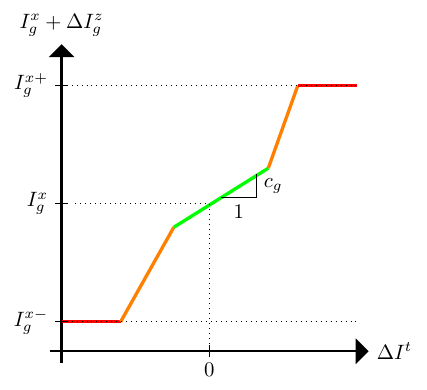}
    \caption{Illustration of the behavior of generator $g$.}
    \label{fig:illustrate_generator}
\end{figure}

As a result, generators will increase their injections at a faster rate once some of the other generators have reached their injection limits. The resulting behavior of the generator injection is visualized in \cref{fig:illustrate_generator}. 
It is difficult to provide an explicit description of the control offset  $ \Delta I_{g}^{z}$ at generator $g$ as a function of the uncertainty offsets $\Delta I _n^y$ for nodes $n\in \mc{N}$ and the generator set points $I_{\bar{g}}^{x}$ for generators $\bar{g}\in \mc{G}$. However, as long as the load distribution control is feasible, i.e., we assume 
\begin{equation}
    \sum_{g\in\mc{G}} -I_{g}^{x+} \le \sum_{n\in \mc{N}} I_n^0 + \Delta I_n^y \le \sum_{g\in\mc{G}} -I_{g}^{x-}\; , \label{eq:unique_condition}
\end{equation}

 this  model for the load distribution implicitly sets the control offsets to unique values. This is an  important characteristic of the model because it means that the control actions $\mb{z}$ do not include any degree of freedom regarding the injections of the network.

\section{Algorithm}
Our algorithm is based on an algorithm \cite{djelassi2021esipglobal} that is suitable for the global solution of so-called existence-constrained semi-infinite optimization problems (ESIPs) with continuous nonlinear functions. 
We show that \cref{eq:flex} can be formulated as an ESIP and specialize the algorithm  based on some of the properties of \cref{eq:flex}.
In the following discussion, we will use $g(\mb{x},\mb{y},\mb{z}):= \max\left\{ g^s(\mb{x},\mb{y},\mb{z},\mb{s}(\mb{x},\mb{y},\mb{z})),g^s(\mb{x},\mb{y}^0,\mb{z}^0,\mb{s}(\mb{x},\mb{y}^0,\mb{z}^0)) \right\}$ to shorten the notation.
\subsection{Reformulation to an existence constrained semi-infinite optimization problem \label{sec:reform}}
Problem \cref{eq:flex} is a slight generalization of an existence-constrained semi-infinite optimization problem (ESIP) which have the form
\begin{equation}
\begin{aligned} \label{eq:ESIP-base}
\inf_{\mb{\check{x}}\in \mc{X}} \; & \check{f}(\mb{\check{x}})\\
\text{s.t. } & \forall \mb{\check{y}} \in \mc{\check{T}}\left[ \exists {\mb{\check{z}}\in \check{\mc{Z}}(\mb{\check{x}},\mb{\check{y}}}) \right. :\\ 
	         & \qquad \left. \check{g}(\mb{\check{x}},\mb{\check{y}},\mb{\check{z}}) \le 0 \right]
\end{aligned}
\end{equation}
and for which a global solution procedure was presented in \cite{djelassi2021esipglobal}. 
 
The problem \cref{eq:flex} is a generalization in the sense that the uncertainty domain $\mc{T}(\mb{x},\delta)$ is parametric in the preventive actions, while $\check{\mc{T}}$ in \cref{eq:ESIP-base} is constant.
However, the following proposition shows that we can use a commonly used relaxation to remove the dependency of the uncertainty domain on the preventive actions. Furthermore, this relaxation is exact in our application. 

\begin{proposition}
Assume that the uncertainty domain is defined by a continuous function $h$ in the form $\mc{T}(\mb{x},\delta) := \{\mb{y}\in \bar{\mc{Y}}| h(\mb{x},\mb{y}) \le \delta\}$ with a given host set for the uncertain values $\bar{\mc{Y}}$. We can exactly reformulate problem \cref{eq:flex} to  a problem with fixed uncertainty set by replacing the constraint
\begin{equation}
\forall \mb{y}\in \mc{T}(\delta,\mb{x}) [\exists \mb{z}\in \mc{Z}(\mb{x},\mb{y}): g(\mb{x},\mb{y},\mb{z})\le 0]
\label{eq:GSIP_feas_constr_orig}
\end{equation}
with
\begin{equation*}
\forall \mb{y}\in \mc{\bar{Y}} [\exists \mb{z}\in \mc{Z}(\mb{x},\mb{y}): \min\left\{\alpha\left(\delta-h(\mb{x},\mb{y})\right),    g(\mb{x},\mb{y},\mb{z})\right\}]    
\end{equation*}
for any arbitrarily chosen but fixed scaling constant $\alpha \in \R^+$.
\begin{proof}
We reformulate the  constraint  for the preventive actions $\mb{x}$ and the flexibility index $\delta$ in problem \cref{eq:flex}, i.e., \cref{eq:GSIP_feas_constr_orig} 
in the following way:
\begin{subequations}
\label{eq:GSIP_REL}
\begin{alignat}{1}
\,&  \forall \mb{y} \in T(\mb{x},\delta) \,[\exists \mb{z} \in \mc{Z}(\mb{x},\mb{y}) :  g(\mb{x},\mb{y},\mb{z})\le 0] \label{eq:GSIP_feas_quant}\\
\Leftrightarrow \, &  \forall \mb{y} \in T(\mb{x},\delta) \,[\mpmin{\mb{z} \in \mc{Z}(\mb{x},\mb{y})} g(\mb{x},\mb{y},\mb{z})\le 0] \label{eq:GSIP_feas_half_quant}\\
\Leftrightarrow \, &  \forall \mb{y} \in \mc{\bar{Y}} \, [ \delta<h(\mb{x},\mb{y})  \lor  \mpmin{\mb{z}\in \mc{Z}(\mb{x},\mb{y})} g(\mb{x},\mb{y},\mb{z}) \le 0 ] \label{eq:GSIP_feas_constr}.\\
\intertext{ 
 We relax the strict inequality and obtain}
\,& \forall \mb{y} \in \mc{\bar{Y}} \,  [  \delta \le h(\mb{x},\mb{y})  \lor   \mpmin{\mb{z}\in \mc{Z}(\mb{x},\mb{y})} g(\mb{x},\mb{y},\mb{z}) \le 0] \label{eq:GSIP_feas_constr_relax} . \\
\intertext{
which in turn is equivalent to
}
        \,& \forall \mb{y} \in \mc{\bar{Y}} \, [ \, \mpmin{\mb{z}\in \mc{Z}(\mb{x},\mb{y})} \mpmin{}\left\{\alpha\left(\delta-h(\mb{x},\mb{y})\right),    g(\mb{x},\mb{y},\mb{z})\right\} \le 0] \\
      \Leftrightarrow \, &0  \ge \mpsup{\mb{y}\in\mc{\bar{Y}}}\mpmin{\mb{z}\in \mc{Z}(\mb{x},\mb{y})}\mpmin{}\left\{\alpha\left(\delta-h(\mb{x},\mb{y})\right),    g(\mb{x},\mb{y},\mb{z})\right\} \label{eq:MLP_flex} \\
      \Leftrightarrow \, &\forall \mb{y} \in \mc{\bar{Y}} [\exists \mb{z} \in \mc{Z}(\mb{x},\mb{y}) :  g(\mb{x},\mb{y},\mb{z})\le 0]  
\end{alignat}      
\end{subequations}
for any arbitrary scaling constant $\alpha \in \R^+$.

It remains to show that the relaxation from \cref{eq:GSIP_feas_constr} to \cref{eq:GSIP_feas_constr_relax} is exact in this specific case, in the sense that the  infimum of \cref{eq:flex}  is the same as the minimum of the relaxation. To see this, consider any solution $(\delta,\mb{x})$ fulfilling the relaxation \cref{eq:GSIP_feas_constr_relax}.  For any $\epsilon>0$ the disturbed solution $(\delta-\epsilon,\mb{x})$ fulfills the original constraint \cref{eq:GSIP_feas_constr} and is thus feasible in the original problem.  If the point $(\delta,\mb{x})$ was a global minimizer in the relaxation, the corresponding objective value is $\delta$. The objective value of the disturbed solution is $\delta-\epsilon$. As $\epsilon$ goes to zero, the relaxation only introduces an arbitrarily small error in the objective function value obtained. 
\end{proof}
\end{proposition}

Note that the set of uncertainties is now $\bar{\mc{Y}}$ and does no longer vary with the preventive actions $\mb{x}$. Consequently, we can use existing algorithms for ESIPs.

\subsection{Algorithmic steps}

We employ an adaptive discretization approach. The fundamental idea is to iteratively add worst-case scenarios for the values of the uncertain variables and dates back to \cite{blankenship1976infinitely}.

In \cite{djelassi2021esipglobal}, we detailed how to extend our previous algorithm from \cite{djelassi2017hybrid} for standard semi-infinite optimization problems to ESIPs. The extended algorithm used there is a modification of our restriction of right-hand side algorithm for semi-infinite optimization problems \cite{mitsos2011global}.

 Instead of the modified algorithm, we use the restriction of the right-hand side algorithm as the basis here, because it is easier to parallelize.  We believe it is necessary to use parallelization and  focus on wall run time, because in application, we might start the calculation based on current information and  we expect that the uncertainty compared to the nominal state increase as time passes.  

The algorithm uses two subproblems, the upper-level and the lower-level problem.
In the upper-level problem, we solve 
\begin{equation}
\begin{aligned} \label{eq:ULP}
\smashoperator[r]{\min_{\delta\in \R^+,\mb{x}\in \mc{X},\mb{z}^1,...,\mb{z}^{|\mc{D}|}}} \; &-\delta \\
\text{s.t. } & \, \forall  d \in \mc{D}: \\
&   \mpmin{}\left\{\alpha\left(\delta-h(\mb{x},\mb{y}^d)\right),    g(\mb{x},\mb{y}^d,\mb{z}^d)\right\}\le -\epsilon_R\\
& \, \forall  d \in \mc{D}: \mb{z}^d\in \mc{Z}(\mb{x},\mb{y}^d)
\end{aligned}
\end{equation}
for a given discretization $\mc{Y}^{disc}=\{\mb{y}^d , d\in \mc{D}\}$ with index set $\mc{D}$ and restriction $\epsilon_R>0$.
We run two procedures in parallel. 
In the lower-bounding procedure, we  set the restriction $\epsilon_R$ to zero which means that every objective value of  \cref{eq:ULP} is an update to the lower bound  $-\delta^{LB}$   for \cref{eq:flex}.
We also run an upper-bounding procedure that uses positive values of $\epsilon_R$, starting from $\epsilon_R^0$, in order to find a feasible solution (and thus an update to the  upper bound $-\delta^{UB}$).
Using the scheme to iteratively reduce $\epsilon_R$ (by a constant factor $\frac{1}{r_R}$) introduced in \cite{mitsos2011global}, a feasible and $\epsilon$-optimal solution can be found in finitely many steps.

In both procedures, we also solve a worst-case generation problem, which is represented by the right-hand-side of \cref{eq:MLP_flex}.
In the  worst-case generation  problem, we solve for the worst-case uncertainty values for given values of $\delta$ and $\mb{x}$.  
If the  worst-case generation  problem has a positive maximal objective value, it provides the next addition to the discretization $\mc{Y}^{disc}$. Otherwise, it confirms the feasibility of  the pair $(\mb{x},\delta)$.

The  worst-case generation  problem \cref{eq:MLP_flex} is a maxmin problem and is solved via a discretization approach for minmax problems \cite{falk1977nonconvex}. This solution approach for the  worst-case generation  problem is similar to the procedure outlined above. However, we do not require an analog to the  upper-bounding procedure, since the objective value of the inner minimization problem already provides the required bound \cite{falk1977nonconvex}.

\subsection{Specialization \label{sec:algorithmic_specializations}}
We specialize the algorithm with three changes. The first two aim to exploit the special case where $\mc{T}(\mb{x},\delta) =  [\mb{y}^0-\mb{\Delta}^-\delta, \mb{y}^0+\mb{\Delta}^+\delta]\subseteq\mc{\bar{Y}}$, which is a scaled hyperbox around a nominal point $\mb{y}^0$ with fixed scaling factor vectors $\mb{\Delta}^+$ and $\mb{\Delta}^-$. The other provides an additional way of finding an upper bound.
\subsubsection{Dropping  redundant discretization points}
When the parametric uncertainty region does not depend on the preventive actions $\mb{x}$, i.e., $\mc{T}(\mb{x},\delta) =\tilde{\mc{T}}(\delta):=\{\mb{y}\in \bar{Y}| \tilde{h}(\mb{y}) \ge \delta\}$, we  can reduce the size of the optimization problem \cref{eq:ULP}, by  removing points $\mb{y}^d$ with $\tilde{h}(\mb{y}^d)\ge \delta^{LB}$ for the following iterations. Indeed, with the information $\delta\le \delta^{LB}$, these discretization points are redundant. Note that each discretization point increases the problem size not only due to the  additional constraints but also due to the additional decision variables $\mb{z}^d$ added to the problem.

\subsubsection{Transformation-based reformulation\label{sec:transformation_reformulation}}
\begin{figure*}[b]\centering
\subfloat[Behavior of the algorithm without the transformation.  The original box $\mc{T}(\delta^0)$ (thick black line) is reduced so that the discretization point $\mb{y}^1$ is outside of the new box $\mc{T}(\delta^1)$ (dotted line). \label{fig:no_proj}] {\includegraphics[width=0.485\textwidth]{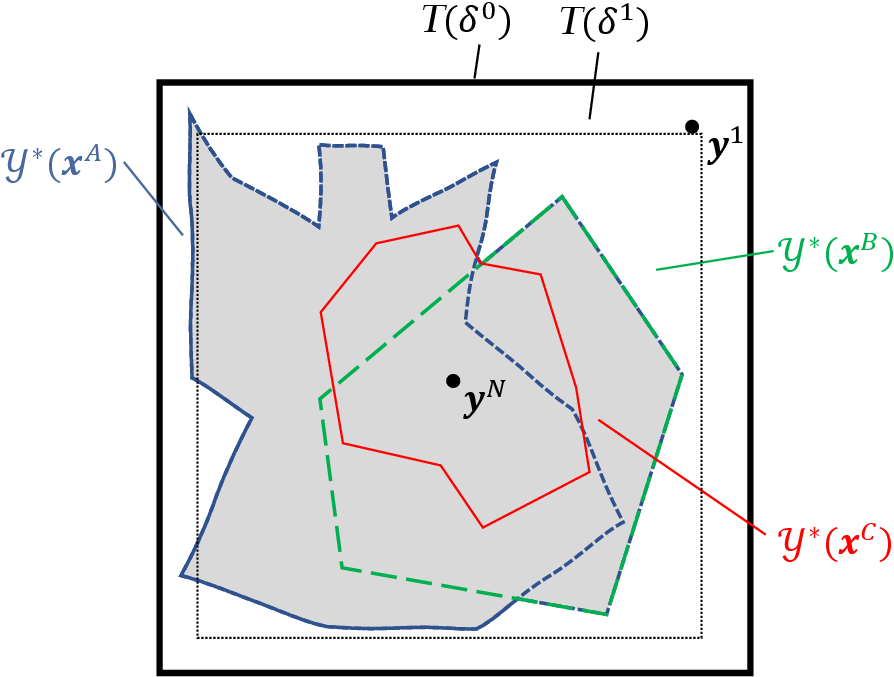}}\quad
\subfloat[Behavior of the algorithm with the transformation. The discretization point $\mb{y}^1$ is projected on the boundary of the box. The box must be  reduced until the transformed point $\mb{y}^{1'}$ can be handled with one of the choices for the preventive actions (gray area). The resulting box $\mc{T}(\delta^1)$ (finely dotted line) is smaller.\label{fig:proj}]{\includegraphics[width=0.485\textwidth]{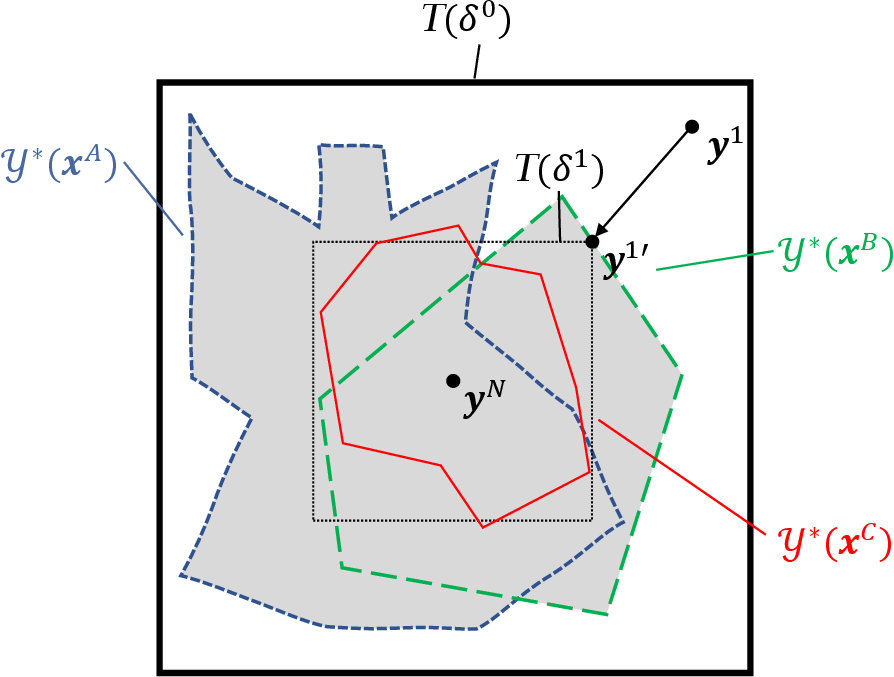}}
\caption{Illustration of the behavior of the algorithm with and without transformation in a fictitious example. $\mc{Y}^*(x)$ denotes the set of all uncertainty values $\mb{y}$ that can be handled with the preventive actions $\mb{x}$. None of the preventive actions $\mb{x}^A,\mb{x}^B,\mb{x}^C$ can handle the discretization point $\mb{y}^1$.}
\label{fig:compare}\end{figure*}

As mentioned above, discretization points become redundant when the lower bound satisfies $\tilde{h}(\mb{y}^d)\ge \delta^{LB}$. 
Similar behavior occurs during the solution of \cref{eq:ULP}.
 Because we maximize $\delta$, its value will only fall below $h(\mb{y}^d)$ if there are no preventive actions $\mb{x}$ and controls $\mb{z}^d$ with $g(\mb{x},\mb{y}^d,\mb{z}^d)\le0$. 
 However,  if $\delta$ falls below $ h(\mb{y}^d)$, the discretization point has no further influence. 
  In the specific case of a scaled hyperbox, i.e., $\tilde{T}(\delta) = [\mb{y}^0-\mb{\Delta}^-\delta, \mb{y}^0+\mb{\Delta}^+\delta]\subseteq\mc{\bar{Y}}$, this would mean that a discretization point $\mb{y}^d$ will be ignored if the hyperbox is reduced enough to exclude it. This is shown in \cref{fig:no_proj}.
However, we can still utilize the discretization point by using a transformation. 
Specifically, we choose to transform it so that it traces a path from the border of the hyperbox instead of being ignored by using 

\begin{equation}
\begin{aligned} \label{eq:ULP_hyperbox_mid}
\smashoperator[r]{\min_{\delta\in \R^+,\mb{x}\in \mc{X},\mb{z}^1,...,\mb{z}^{|\mc{D}|}}} \; &-\delta \\
\text{s.t. } & \, \forall  d \in \mc{D} : \\
&g(\mb{x},\mb{y}^0+(\mb{y}^d-\mb{y}^0) \frac{\min\{h(y^d),\delta\}}{h(y^d)},\mb{z}^d) \le -\epsilon_R.\\
& \, \forall  d \in \mc{D}: \mb{z}^d\in \mc{Z}(\mb{x},\mb{y}^d)
\end{aligned}
\end{equation}
instead of \cref{eq:ULP}. This is visualized in \cref{fig:proj}.
When using this transformation, we can choose to drop the discretization points as mentioned above, but it can also be beneficial to keep them.
Similar transformations are known in the context of generalized semi-infinite programming \cite{still1999generalized}.  A similar idea has also been discussed in \cite{zhao2022novel} in the context of flexibility index evaluation without discretization-based methods.
The expected benefit of using this reformulation is reducing the total number of necessary iterations by  reducing $\delta$ further in a single iteration.

\subsubsection{Upper-bounding heuristic \label{sec:upper_bounding_heuristic}}
In the procedure described, an upper bound is found when the upper or lower-bounding procedure produces a pair ($\mb{x},\delta$) that is feasible, i.e., fulfills \cref{eq:MLP_flex}. 
Note that  \cref{eq:ULP} can return a pair  of preventive actions $\mb{x}$ and  flexibility index $\delta$  which is infeasible but would be feasible for a reduced $\delta$.
Also note that we can rewrite the feasible region of \cref{eq:MLP_flex} as
\begin{equation}
\begin{split}
\Leftrightarrow \,& \delta\le  \mpinf{\mb{y} \in \mc{\bar{Y}}} h(\mb{x},\mb{y}) \\
 \, & \quad  \quad \text{s.t.} \mpmin{\mb{z}\in \mc{Z}(\mb{x},\mb{y})} g(\mb{x},\mb{y},\mb{z})>0.
\end{split}\label{eq:MLP_heur}
\end{equation}
Solving the problem on the right-hand side of the inequality of \cref{eq:MLP_heur}, we can find  $-\delta_{wc}(\mb{x})$, a tight upper-bound on our objective value in \cref{eq:flex}, $-\delta$, for fixed preventive actions $\mb{x}$. A similar problem for finding the flexibility index for fixed preventive actions $\mb{x}$ in a more restricted setting appears in \cite{grossmann1987active}. In our setting, finding the flexibility index for fixed preventive actions $\mb{x}$, i.e, solving \cref{eq:MLP_heur}, is only used as a subproblem, but has useful applications on its own. Variants of this problem are discussed in the power flow literature, see for example \cite{zhao2015unified,capitanescu2021power}.
When the strict inequality is relaxed, this problem is a semi-infinite optimization problem. In our test cases, we again use the restriction of the right-hand side algorithm from \cite{mitsos2015globalGSIP} to solve this relaxed problem.
Since we relax  problem  \cref{eq:MLP_heur}, we do not necessarily obtain the tight upper-bound $-\delta_{wc}(\mb{x})$, but a pessimistic value $-\delta_{wc}^{relax}(\mb{x})\ge -\delta_{wc}(\mb{x})$. 
This procedure is started in parallel when new preventive actions $\mb{x}$ are encountered in either the lower or upper-bounding procedure.
\section{Parametric uncertainty region formulations}

A key aspect in choosing the parametrization of the uncertainty region $\mc{T}(\mb{x},\delta)$ is the interpretability of the resulting value of $\delta$.
 To illustrate the potential uses of the formulation \eqref{eq:flex}, we investigate two parametrizations for the uncertainty region in the context of flexibility of power grids. They differ mainly in the choice of the parametric uncertainty region $\mc{T}(\mb{x},\delta)$. In both cases, the  value of $\delta$ has a specific interpretable meaning.

\subsection{Inner approximation of region of manageable uncertain injections with a scaled hyperbox\label{sec:box_param}}
In this case, we are given  the forecasted uncertain injections $\mb{y}^0 $ and scaling vectors $\mb{\Delta}^-$  and $\mb{\Delta}^+$ for the expected positive and negative deviation from that forecast, respectively. We are also given an initial known upper bound for the flexibility $\delta^{UB}$ that defines the host set for the uncertainties $\bar{\mc{Y}}$. We search for the largest flexibility index $\delta$ and decision variables $\mb{x}$ (representing generator set-points), such that all uncertain injections  $\tilde{T}(\delta) = [\mb{y}^0-\mb{\Delta}^-\delta, \mb{y}^0+\mb{\Delta}^+\delta]$ can be handled while ensuring that the power flow $P_e$ over on the critical edges $e\in \mc{E}_c$ does not exceed a specified limit $P_e^{lim}$, i.e., $g$ represents the constraints

\begin{equation} \label{eq:overload_const_multiple}
\mabs{P_e} \le {P_e^{lim}}, \quad \forall e \in \mc{E}_c
\end{equation}

which we can summarize as the scalar constraint
\begin{equation} \label{eq:overload_const}
\mpmax{e \in \mc{E}_c}\frac{\mabs{P_e}}{P_e^{lim}}- 1\le 0 .
\end{equation}

The resulting hyperbox will constitute an inner approximation of the region of all manageable uncertain injections. Its size is maximized by the returned preventive actions $\mb{x}$. A possible application of this parametrization is to search for the maximal uncertainty of node injections relative to the initial or forecasted injections.
\subsection{Maximal available net power transfer capacity\label{sec:max_flow_param}}
We also investigate  the problem of finding \maxflow{} $\, \delta$ from a region $A$ to a region $B$ of the network.
We assume that we are given a range for potential uncertain injections, i.e., $\bar{\mc{Y}}=[\mb{y}^-,\mb{y}^+]$ and forecasted values $\mb{y}^0$.
 For the found power transfer capacity $\delta$, we must guarantee  nominal grid operations, i.e., \cref{eq:overload_const},
 for any realization of the uncertainty that leads to an additional net power transfer within $[0,\delta]$ from region $A$ to a region $B$.
Note that this differs from finding the maximal additional power transfer $\delta^{max}$ that is can be safely achieved under the specified uncertainty, as this would not guarantee that \emph{all} power transfers $\delta^\dagger$ with $0\le\delta^\dagger<\delta^{max}$ can be archived safely.

We formulate this as follows:
We denote with $\mc{N}_{A}$ and $\mc{N}_{B}$ the sets of nodes inside the regions $A$ and $B$ regions and with $\mc{G}_{A}$ and $\mc{G}_{B}$ the generators within these regions, respectively.
The  additional net power transfer is the minimum of the additional injections in the $A$ region and the decrease of injections in the $B$ region.
As a result, the function $h$ describing the uncertainty region takes the form
\begin{equation}
	\begin{split}
	h(\mb{x},\mb{y})= \mpmin{}\left\{\sum_{n\in \mc{N}_{A}} I_n-I_n^0 - \sum_{g \in \mc{G}_{A}} I_g^x \right. ,\\ \left.  \sum_{n\in \mc{N}_{B}} -I_n+I_n^0 + \sum_{g \in \mc{G}_{B}} I_g^x   \right\}. \label{eq:net_pow_transfer}
	\end{split}
\end{equation}
Since we only need to handle uncertainty values $\mb{y}$ resulting in positive flows from $A$ to $B$, the constraint $g$ takes the form
\begin{equation}
    \min\{\alpha h(\mb{x},\bm{y}),\mpmax{e \in \mc{E}_c}\frac{\mabs{P_e}}{P_e^{lim}}- 1\}\le 0  \;.
\end{equation}

Note that, as mentioned earlier, the control offsets for generators $g\in \mc{G}$ are uniquely determined by the generator set points $I_{{g}}^{x}$ at all generators $g\in\mc{G}$ (included in the preventive actions $\mb{x}$) and the additional injections  $\Delta I_n^y$ at all nodes $n \in \mc{N}$  (included in the uncertain values $\mb{y}$). As a result, the injections $I_n$ for $n\in\mc{N}$ are  not a function of the control variables $\mb{z}$.
Therefore, the net flow between the regions as defined by \cref{eq:net_pow_transfer} is also independent of the control variables $\mb{z}$.

\section{Numerical Experiments}

We now report results from numerical experiments for both  parametric uncertainty region formulations on a small-scale instance based on the modified 30-bus IEEE system from \cite{djelassi2018hierarchicalgrid} and a medium-scale grid instance generated at RTE.

All computations are run on a system with a Skylake Platinum 8160 processor using 24 cores and 48GB of memory.
To solve the resulting mixed-integer optimization problems, we use Gurobi 10.0.2. Each instance of the MILP solver is allowed to use up to 6 threads. 
For the medium-sized instance, numerical problems were caused
on some runs by integrality violations. Instead of reducing the integrality tolerance from its default
value, we used the \emph{IntegralityFocus} setting in Gurobi.
We use a relative optimality tolerance of $0.05$ for the overall algorithm and $0.025$ for the auxiliary problem \cref{eq:MLP_heur}. The initial restriction $\epsilon^0_R$ is set to $0.05$ for the overall algorithm and $0.005$ for the auxiliary problem. In both cases, we use $r_R=2$ to determine the reduction rate of the restriction $\epsilon_R$.

For the parametrization from  \cref{sec:box_param}, we use the largest number that guarantees \cref{eq:unique_condition} as the initial upper bound $\delta^{UB}$.

\subsection{Motivating example}

\begin{figure}
  \centering
  {\scalebox{1.00}{\input{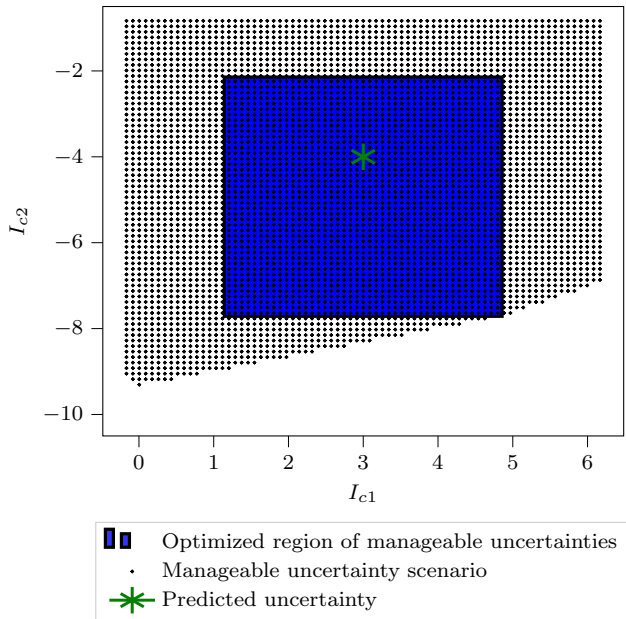}}}
  \caption{Visualization of the result from the motivating example for the final preventive actions. The optimized region of manageable uncertainties is valid and optimal for the chosen parameterization. \label{fig:motivating_res}}
\end{figure}

We solve the motivating example from \cref{sec:motivating_example} with the parameters $h=\SI{1}{\siemens}$, ${P}^{lim}=\SI{5}{\mega \watt}$ $I_{c1}^{pred}=\SI{3}{\mega \watt}$, $I_{c2}^{pred}=\SI{-4}{\mega \watt}$ and $I_{g1}\in [\SI{-7.5}{\mega \watt},\SI{3}{\mega \watt}]$, $I_{g2}\in [\SI{-3}{\mega \watt},\SI{7.5}{\mega \watt}]$. 
Here, we are interested in the uncertainty parameterization as a scaled hyperbox from \cref{sec:box_param}.
Note that the scaling vectors given in the description of the example are $\bm{\Delta}^{+}=[1,1]^T$ and $\mb{\Delta}^{-}=[1,2]^T$.
This problem is solved in approximately 1 second and yields $\delta \in [-1.857,-1.857]$.
The low dimensionality of the problem allows us to visualize the resulting hyperbox for the final values for the generator set-points $I_{g1}^0$ and $I_{g2}^0$ evaluate for  sampled scenarios for the values of uncertain injections if they would lead to a constraint violation.
The results is shown in \cref{fig:motivating_res}.
We see that in this example and for the final values of the preventive actions (here the generator set-points), the actual region of manageable uncertainties is a convex polyhedron, but we want to emphasize this does not hold in general.

The simple scaled hyperbox parameterization allows for an easily interpretable result and allows for faster computation compared to other parameterizations, such as a box where all dimensions are scaled individually.
However, the visualization also shows that the result is conservative, i.e., there are much more manageable uncertainties than the ones that are guaranteed to be manageable by the result of the computation.

In cases where a better approximation of the region of manageable uncertainties is desired,  a possible improvement could be to \emph{optimize} with a simple parametrization such as the scaled hyperbox and then \emph{evaluate} for fixed preventive action for another, less conservative parameterization.

\subsection{Small-scale instance \label{sec:small_grid}}
\begin{figure}[tb]
  \centering
\includegraphics[width=0.70\textwidth]{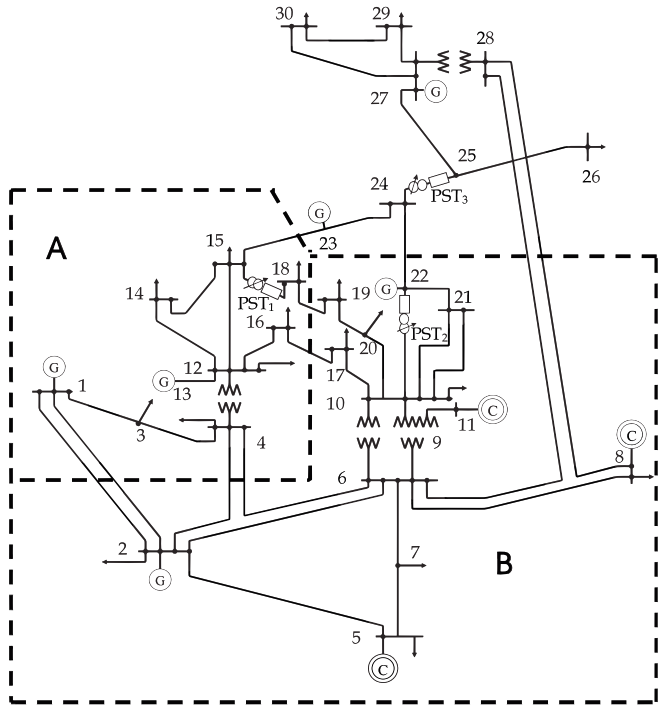}
  \caption{ Visualization of the small grid instance with two regions marked. \label{fig:illustrate_small_grid}}
\end{figure}
As a small-scale instance, we use the grid instance derived in \cite{djelassi2018hierarchicalgrid} based on the 30-bus IEEE system. \cref{fig:illustrate_small_grid} shows the topology of the grid. The main modifications compared to the 30-bus IEEE system are the introduction of three PSTs and the splitting of 7 nodes that are then subject to the bus merging.
In \cite{djelassi2018hierarchicalgrid}, it was shown that if injections at non-generator nodes are considered to be
uncertain within a range of $\pm30\%$ of the nominal injections, all uncertainty realizations are manageable. To avoid the situation where the algorithm finishes in the first iteration for the parametric uncertainty formulation in \cref{sec:max_flow_param}, we increased the range to $\pm45\%$ and increased the upper generator bounds by a factor of $1.13$ to avoid a violation of \cref{eq:unique_condition}. \Cref{tab:summary_small_grid} summarizes key properties of this grid instance.
\begin{table}[htb]
  \begin{center}
    \caption{Properties of the small grid instance}
    \label{tab:summary_small_grid}
    \begin{filecontents*}{./small_case_meta_grid.txt}
#   supernodes  nodes   generators  edges   overloads   subedges    shifters     regions
    29           37      6           41      41          7           3		   3	
\end{filecontents*}
\pgfplotstabletypeset[
      alias/critical edges/.initial=overloads,
      alias/bus merges/.initial=subedges,
      multicolumn names, ignore chars={\#},
      columns={nodes, generators,edges,shifters,critical edges,bus merges},
      col sep=space, every head row/.style={before row={\toprule}},   
        every last row/.style={after row=\bottomrule}, every column/.style={sci,sci zerofill, precision=2,int detect}
    ]{./small_case_meta_grid.txt} \end{center}
\end{table}
\subsubsection{Influence of $\alpha$}
We investigate the influence of the value of parameter $\alpha$ in \cref{eq:MLP_flex} on the solution time.

The selected worst-case uncertainty values $\mb{y}^*$, i.e., the solutions to the maxmin problem \cref{eq:MLP_flex} that are used as discretization points, are selected according to two objectives.
The parameter $\alpha$ weights the two objectives $\delta-h(\mb{x},\mb{y}^*)$ and $\mpmin{\mb{z}\in \mc{Z}(\mb{x},\mb{y}^*)} g(\mb{x},\mb{y}^*,\mb{z})$. The first can be understood to measure how far the discretization point $\mb{y}^*$ is inside the parametric uncertainty region $\mc{T}(\delta,\mb{x})$, while the second can be understood to measure the maximum violation considering the available control variables.
For very small values of $\alpha$, the overall objective will be determined more by the first objective. 

The value of $g$ in \cref{eq:overload_const} is already normalized. In the following, we will choose $\alpha=\alpha^{'}/\delta^{norm}$ to achieve a similar effect. For the parametrization of \cref{sec:box_param}, we use $\delta^{norm}=1$ and  for the parametrization of \cref{sec:max_flow_param} we use choose $\delta^{norm}$ as the optimistic value for $\delta$, i.e., the value obtained when the uncertain injections would also become decision variables.

\begin{table*}
\centering
\begin{tabular}{lllll}
\hline
\multicolumn{1}{l|}{$\alpha$} & all & \thead{without \\ auxiliary problem} & \thead{without \\ dropping points} & \thead{without \\ transformation} \\ \hline
\multicolumn{1}{l|}{0.1} & TIME\_OUT   & TIME\_OUT & TIME\_OUT   & TIME\_OUT    \\
\multicolumn{1}{l|}{0.5} & 50-53          & 2850-2855      & 357-362         & \textbf{34-35}  \\
\multicolumn{1}{l|}{1}   & \textbf{77-80} & 136-138       & 195-199         & 1181-1194         \\
\multicolumn{1}{l|}{2}   & 599-608         & 580-582       & 956-960         & \textbf{122-126} \\
\multicolumn{1}{l|}{10}  & 163-165         & 149-150       & \textbf{25-27} & TIME\_OUT    \\ \hline
\end{tabular}

\caption{Range of wall run times in seconds for the small grid instance with the parametrization from \cref{sec:box_param} for different values of $\alpha^{'}$ . TIME\_OUT denotes that the run-time limit of 3600s was exceeded. \label{tab:alpha_box} The best time is given in bold font. Not using the one of  proposed specializations  leads to a significant increase in computational time for at least one value of $\alpha^{'}$. Run times vary strongly with $\alpha^{'}$.}
\end{table*}

\begin{table*}
\centering
\begin{tabular}{lll}
\hline
\multicolumn{1}{l|}{$\alpha$} & all & without auxiliary problem   \\ \hline
\multicolumn{1}{l|}{0.05} &\bf{1228-1234}   & 1714-1838     \\
\multicolumn{1}{l|}{0.1} & \bf{587-618}  & 1106-1186     \\
\multicolumn{1}{l|}{0.5} & \bf{819-825}        & 2145-2247    \\
\multicolumn{1}{l|}{1}   & TIME\_OUT & TIME\_OUT  \\
\multicolumn{1}{l|}{2}   & TIME\_OUT   & TIME\_OUT     \\
\multicolumn{1}{l|}{10}  &TIME\_OUT   & TIME\_OUT     \\ \hline
\end{tabular}

\caption{Range of wall run times in seconds for the small grid instance with the parametrization from \cref{sec:max_flow_param} for different values of $\alpha^{'}$ . TIME\_OUT denotes that the run-time limit of 3600s was exceeded. \label{tab:alpha_flow} The best time is given in bold font. Not using the the  proposed specializations of solving the auxiliary problem consistently leads to a significant increase in computational time for at least one value of $\alpha'$. Run times vary strongly with $\alpha^{'}$.}
\end{table*}

\Cref{tab:alpha_box} and \cref{tab:alpha_flow} shows the range of wall times on the small grid instance for different values of $\alpha'$, for the parametrization from \cref{sec:box_param} and \cref{sec:max_flow_param}, respectively, for different configurations of the algorithm. The time range given is taken over three repetitions.

We  observe that the run time is strongly affected by the value of $\alpha'$, where the run time can more than double from one choice of $\alpha'$ to the next. This phenomenon is independent of the application of the algorithmic specializations from \cref{sec:algorithmic_specializations}. 

For the box parameterization there does not seem to be a common good choice of $\alpha'$ between the different variants of the proposed algorithm. There is also no configuration that consistently outperforms the others. However, the variant using all the proposed specialization seems to be the most consistent over the different values of $\alpha'$  with the others having one value of $\alpha'$ where they took much longer than the other configurations. 
Not using the transformation specialization was either very beneficial with it being the fastest approach or lead to extremely increased computation times. Clearly, forgoing the transformation for simpler subproblems can be beneficial if the number of iterations stays low.
We were surprised to find that for some values of $\alpha'$, not solving the auxiliary problem was faster, even though its computation was run in parallel. However, it turns out that in those instances, the auxiliary problem did simply not matter in terms of needed iterations and that the additional time was spent waiting for the abortion of the solution process of the auxiliary problem. This could potentially be prevented with a more forceful termination procedure. For the parameterization from \cref{sec:max_flow_param}, where the run times were generally longer, using the auxiliary problem was always beneficial.

Overall, the algorithmic framework exhibits a significant inherent variation in run time with variations in $\alpha$ and similarly $P^{lim}$. This is due to the complex alteration of solution paths caused by the strong dependence of one subproblem solution on the sequence of earlier subproblem solutions.

Still, further research to explore whether an optimal static value of $\alpha$ for a given problem can be determined in advance or if adapting $\alpha$ during the algorithm execution would be beneficial.

In the absence of a systematically good choice, we choose $\alpha'=0.5$ for the the small grid instance.

\subsubsection{Results}
First, we use the parametrization  from  \cref{sec:max_flow_param} in formulation \cref{eq:flex} to compute the \maxflow{} from part A to B of the network for the specified uncertainty ranges of $\pm 45\%$ of the nominal injections. The regions A and B are shown in the \cref{fig:illustrate_small_grid}. Using $\alpha'=0.5$,
we terminate with $\delta\in[\SI{0.590}{\mega \watt},\SI{0.614}{\mega \watt}]$ after approximately 780 seconds using 18 iterations of the lower and 16 iterations of the upper-bounding procedure, respectively.

Second, we use the parametrization from \cref{sec:box_param} with $\bm{\Delta}^{+}$ and $\mb{\Delta}^{-}$ chosen such that $\delta=1$ represents the aforementioned $\pm 45\%$ deviation of the nominal injections. Using $\alpha'=0.5$ yields a  flexibility index $\delta\in [0.925,0.954]$ in 12 iterations of the lower and the upper-bounding procedure, respectively, and approximately 63 seconds. We can thus guarantee that uncertainties at around $41\%$ of the nominal injections can be handled. Repeating the computation of  the \maxflow{} for this reduced amount of uncertainty ($41\%$ instead of $45\%$), we obtain a capacity of $\SI{0.798}{\mega \watt}$. In this case, optimizing with a lower amount of uncertainty allows for a larger guaranteed available power transfer.
\subsection{Medium-scale instance}
\begin{table}[ht]
  \begin{center}
    \caption{Properties of the medium grid instance}
    \label{tab:summary_medium_grid}
    \begin{filecontents*}{./medium_case_meta_grid.txt}
#   supernodes  nodes   generators  edges   overloads   subedges    shifters    regions
    10          6800    32          9398    263           0           2           10 
\end{filecontents*}
\pgfplotstabletypeset[
     alias/critical edges/.initial=overloads,
      alias/bus merges/.initial=subedges,
      multicolumn names, ignore chars={\#},
      columns={nodes, generators,edges,shifters,critical edges,bus merges},
      col sep=space, every head row/.style={before row={\toprule}},   
        every last row/.style={after row=\bottomrule}, every column/.style={sci,sci zerofill, precision=2,int detect}
    ]{./medium_case_meta_grid.txt} \end{center}
\end{table}
The medium-scale models the power transmission grid over France and neighboring countries. Case files are accessible at \cite{MediumScale:828516}.
 In our case study, we simulate  a situation where uncertain injections in Spain are to be handled by generators in France.  All uncertain injections except for one are situated in Spain and all active generators ($g\in G$ with $c_g\neq 0$)  are situated in France and we  want to maximize the \maxflow{} from France to Spain. 
The grid instance also contains nodes situated in other countries, but these do not contain active generators or uncertain injections. An overview over the properties of the  instance is given in \cref{tab:summary_medium_grid}.
\subsubsection{Results}
Using $\alpha'=10$ and the pa\-ra\-me\-tri\-za\-tion from \cref{sec:max_flow_param}, we obtain a maximal flow of approximately $\delta\in[\SI{4530}{\mega \watt},\SI{4760}{\mega \watt}]$ in approximately 300 seconds in 3 iterations of the lower and 2 of the upper-bounding procedure. Here we see that if only a small number of iterations is required until convergence, even large instances can be solved quickly.

For the parametrization from \cref{sec:box_param} and using $\alpha^{'}=10$, we obtain $\delta\in[0.546,0.573]$ in approximately 83 minutes with 17 iterations of the lower and 18 iterations of the upper-bounding procedure. As expected, the computational time required for the medium-scale  instance is significantly larger than for the small-scale instance if a similar number of iterations is needed. Still, because we can utilize sophisticated MILP solvers for the subproblems, the problem can be solved in a reasonable amount of time.  Repeating the computation of  the \maxflow{} for the reduced amount of uncertainty, we obtain a \maxflow{} of $\SI{2600}{\mega \watt}$. In this case, directly optimizing for it allowed for a larger \maxflow{}.

\section{Conclusion}
We formulated the maximization of the robustness of power grid operation as a hierarchical optimization problem, following the flexibility concept from process systems engineering \cite{grossmann1983optimization}. The formulation is rigorous in that it treats the uncertainty explicitly, but the resulting problem requires special solution methods.

For the solution, we specialized a discretization algorithm for existence constrained semi-infinite optimization problems based on \cite{mitsos2011global}. This allows us to include integer variables in the control variables and the uncertainty, for example, to describe piecewise linear functions. In the context of our grid model, which is based on the DC flow approximation, it also requires only the solution of mixed-integer linear optimization problems, which allows us to utilize commercial MILP solvers. As a result, numerical experiments on a medium-sized grid instance were shown to be solvable in reasonable computation times. 

However, we have seen that the solution times can be high even for small grid instances and vary drastically with different choices of the scaling parameter $\alpha$. For larger instances, systematically testing different values is prohibited by the already rather high computational cost. Further research is needed if robust choices for the algorithmic parameters can be found prior to solving the problem.

To illustrate the usefulness of the flexibility framework, we looked at two different interpretations of flexibility for power systems.
The problem of finding a maximal available power transfer capacity between two
parts of the network and finding the maximal manageable uncertainty when parametrized as a hyperbox. We showed results on a small and medium-scale power grid instance. 

The presented methodology can be adapted to  other parametrizations for the uncertainty region, as long as the defining function $h$ is continuous. 
In the results of the medium-sized instance, we have illustrated that this can be vital: simply optimizing the flexibility in terms of a scaled hyperbox allows us to also guarantee a flexibility in terms of \maxflow{} for the reduced amount of uncertainty, but despite reducing the host set of the uncertainty, this yielded worse flexibility than optimizing it directly on the original host set of values for the uncertainties.

In the present work, we only considered the DC flow approximation. We rigorously consider the worst-case under this assumption, but if the real grid exhibits significant nonlinear behavior, this will likely be missed in our model. 
Another more accurate model could, at least in principle, be used, as long as the grid state is uniquely determined by the preventive actions $\mb{x}$, the uncertain values $\mb{y}$ and the controls $\mb{z}$. Although the computational effort can be expected to increase, the application to small grid instances seems tractable. Beyond that, future improvements could enable medium-sized instances for nonlinear models.

In conclusion, we believe that the presented methodology constitutes a valuable approach for increasing the robustness and reliability of power grids, today for DC flow approximation and in the future for nonlinear models.

\section{Statements and Declarations}
We declare the following:
\subsection{Funding}
This research is funded by R\'eseau de transport d'electricit\'e (RTE, France) through the project “Hierarchical Optimization for Worst-Case Analysis of Power Grids”.
\subsection{Competing Interests}
We declare that we have no competing interests.
\subsection{Data availability}
Instance data for the medium size test-case is accessible at \cite{MediumScale:828516}. Further instance data to support the findings of this study are available from the corresponding author upon reasonable request.
\subsection{Ethical Approval}
Ethical approval is not applicable for the content of this manuscript.

\bibliographystyle{spbasic} 
\bibliography{008_X+bib}

\end{document}